\documentclass{gtart_h}
\def\ifplaintex{\expandafter\ifx\csname documentclass\endcsname\relax}

\def\gtp{{\mathsurround=0pt\it $\cal G\mskip-2mu$eometry \&\ 
$\cal T\!\!$opology $\cal P\!$ublications}}  % GT publications

\def\recd{{\small Received:\qua\receiveddate\ifx\reviseddate\relax
\else\qquad Revised:\qua\reviseddate\fi\par}} 

\def\lognumber#1{\def\thelognumber{#1}}
\def\volumenumber#1{\def\thevolumenumber{#1}}
\def\volumeyear#1{\def\thevolumeyear{#1}}
\def\papernumber#1{\def\thepapernumber{#1}}
\def\pagenumbers#1#2{\def\startpage{#1}\def\finishpage{#2}}
\def\published#1{\def\publishdate{#1}}

\def\received#1{\def\receiveddate{#1}}
\def\revised#1{\def\reviseddate{#1}}
\def\accepted#1{\def\accepteddate{#1}}

\let\\\par\let\thelognumber\relax\let\thevolumenumber\relax
\let\thepapernumber\relax\let\thevolumeyear\relax\let\startpage\relax
\let\finishpage\relax\let\publishdate\relax\let\receiveddate\relax
\let\reviseddate\relax\let\accepteddate\relax\let\theasciititle\relax
\let\theasciiauthors\relax
\let\theasciiabstract\relax

\let\theasciiemail\relax

\ifplaintex
\font\logobig=cmssbx10 scaled 3836
\font\logomed=cmssbx10 scaled 2557
\else
\font\logobig=cmssbx10 scaled 4200
\font\logomed=cmssbx10 scaled 2800
\fi

\long\def\makeagttitle{   %%% start of definition of \makeagttitle
\count0=\startpage
\agt\hfill      %   Journal title (top left) 
\hbox to 45truept{\vbox to 0pt{\vglue -13truept{\logomed A\kern -.37em{\logobig 
T}\kern -.38em G}\vss}\hss}
\break
{\small Volume \thevolumenumber\ (\thevolumeyear)
\startpage--\finishpage\nl
Published: \publishdate}

\vglue .25truein

{\parskip=0pt\leftskip 0pt plus
1fil\def\\{\par\smallskip}{\Large\bf\thetitle}\par\medskip} \vglue
0.05truein

{\parskip=0pt\leftskip 0pt plus 1fil\def\\{\par}{\sc\theauthors}
\par\medskip}%
 
\vglue 0.03truein 

{\small\leftskip 25truept\rightskip 25truept{\bf Abstract}\stdspace\theabstract

{\bf AMS Classification}\stdspace\theprimaryclass
\ifx\thesecondaryclass\relax\else; \thesecondaryclass\fi\par
{\bf Keywords}\stdspace \thekeywords\par}\vglue 7truept

}   %%%% end of definition of \makeagttitle

\ifplaintex
\font\phead=cmsl9 scaled 950
\font\pnum=cmbx10 scaled 913
\font\pfoot=cmsl9 scaled 950
\headline{\vbox to 0pt{\vskip -4.5mm\line{\small\phead\ifnum
\count0=\startpage ISSN 1472-2739 (on-line) 1472-2747 (printed)
\hfill {\pnum\folio}\else\ifodd\count0\def\\{ }% 
\ifx\theshorttitle\relax\thetitle\else\theshorttitle\fi\hfill{\pnum\folio}
\else\def\\{ and }{\pnum\folio}\hfill\ifx\theshortauthors\relax\theauthors
\else\theshortauthors\fi\fi\fi}\vss}}
\footline{\vbox to 0pt{\vglue 0mm\line{\small\pfoot\ifnum\count0=\startpage
\copyright\ \gtp\hfill\else
\agt, Volume \thevolumenumber\ (\thevolumeyear)\hfill\fi}\vss}}
\else
\font=cmsl9 scaled 1050
\font\lnum=cmbx10 
\font=cmsl9 scaled 1050
\makeatletter
\def\@oddhead{{\small\lhead\ifnum\count0=\startpage ISSN 1472-2739 
(on-line) 1472-2747 (printed)\hfill {\lnum\number\count0}\else\ifodd\count0
\def\\{ }\ifx\theshorttitle\relax \thetitle \else\theshorttitle\fi\hfill
{\lnum\number\count0}\else\def\\{ and }{\lnum\number\count0}
\hfill\ifx\theshortauthors\relax 
\theauthors\else\theshortauthors\fi\fi\fi}}\def\@evenhead{\@oddhead}
\def\@oddfoot{\small\lfoot\ifnum\count0=\startpage\copyright\ \gtp\hfill\else
\agt, Volume \thevolumenumber\ (\thevolumeyear)\hfill\fi}
\def\@evenfoot{\@oddfoot}
\makeatother
\fi
\let\maketitlepage\makeagttitle
\let\makeshorttitle\maketitlepage
\let\maketitle\maketitlepage

\newwrite\gtoutfile
\long\gdef\makeheadfile{  %%% start of definition of \makeheadfile
{\def\\{, }\def\s{ }
\immediate\openout\gtoutfile head.xxx
\immediate\write\gtoutfile{Proxy-for: \ifx\theasciiauthors\relax
\theauthors\else\theasciiauthors\fi\s<\ifx\theasciiemail\relax\theemail\else\theasciiemail\fi>}
\immediate\write\gtoutfile{\noexpand\\}
\immediate\write\gtoutfile{Authors: \ifx\theasciiauthors\relax
\theauthors\else\theasciiauthors\fi}
{\def\\{ }\immediate\write\gtoutfile{Title: \ifx\theasciititle\relax
\thetitle\else\theasciititle\fi}}
\immediate\write\gtoutfile{Subj-class: GT or SG, GR etc}
\immediate\write\gtoutfile{MSC-class: \theprimaryclass\ifx\thesecondaryclass\relax\else, \thesecondaryclass\fi}
\immediate\write\gtoutfile{Journal-ref: Algebr. Geom. Topol. \thevolumenumber\s
(\thevolumeyear) \startpage-\finishpage}
\immediate\write\gtoutfile{Comments: Published by Algebraic and
Geometric Topology at}
\immediate\write\gtoutfile{\s\s\s  http://www.maths.warwick.ac.uk/agt/AGTVol\thevolumenumber/agt-\thevolumenumber-\thepapernumber.abs.html}
\immediate\write\gtoutfile{\noexpand\\}
\immediate\write\gtoutfile{}
\ifx\theasciiabstract\relax
\immediate\write\gtoutfile{\theabstract}\else
\immediate\write\gtoutfile{\theasciiabstract}\fi
\immediate\write\gtoutfile{}
\immediate\write\gtoutfile{\noexpand\\}
\immediate\write\gtoutfile{}
\immediate\closeout\gtoutfile}}  %%% end of definition of \makeheadfile

\def\maketitlepage{\makeagttitle\makeheadfile}
\let\makeshorttitle\maketitlepage
\let\maketitle\maketitlepage

\lognumber{64}
\volumenumber{5}
\volumeyear{2005}
\papernumber{64}
\pagenumbers{1585}{1587}
\received{15 Febrary 2005} 
\revised{22 October 2005}  
\accepted{24 October 2005}
\published{24 November 2005}

\begin{document}
\title{The topological Hawaiian earring group does not\\embed in the 
inverse limit of free groups}
\shorttitle{The topological Hawaiian earring group}
\author{Paul Fabel}

\address{Drawer MA, Department of Mathematics and 
Statistics\\Mississippi State University, Mississippi State, MS 39762, USA}

\email{fabel@ra.msstate.edu}

\urladdr{http://www2.msstate.edu/~fabel/}

\begin{abstract}
Endowed with natural topologies, the fundamental group of the Hawaiian
earring continuously injects into the inverse limit of free
groups. This note shows the injection fails to have a continuous
inverse. Such a phenomenon was unexpected and appears to contradict
results of another author.
\end{abstract}


\maketitle

\section{Introduction}

Quite generally the based fundamental group $\pi _{1}(X,p)$ of a space $X$
becomes a topological group whose topology is invariant under the homotopy
type of the underlying space $X$ (Corollary 3.4 \cite{Biss}). In the context
of spaces complicated on the small scale the utility of this invariant is
emerging. For example topological $\pi _{1}$ has the potential to distinguish
spaces when the algebraic homotopy groups fail to do so \cite{Fab3}.
Unfortunately even in the simplest cases the topological properties of $%
\pi _{1}(X,p)$ can be challenging to understand.

Consider the familiar Hawaiian earring $X=\cup _{n=1}^{\infty }S_{n}$, (the
union of a null sequence of simple closed curves $S_{n}$ joined at a common
point) and the canonical homomorphism $\phi :\pi _{1}(X)\rightarrow
\lim_{\leftarrow }\pi _{1}(\cup _{i=1}^{n}S_{i})$ $.$

The paper \cite[page 370]{Biss} seems to claim that $\phi $ is also a homeomorphism
onto its image (``$\psi ^{-1}$ is surely continuous as well\dots'').
The intent of this note is to show that such a claim is false. The
monomorphism $\phi $ is not a homeomorphism onto its image, and thus $\phi $
fails to be a topological embedding (Theorem \ref{main}). To prove this we
consider the sequence $[(y_{1}\ast y_{n}\ast y_{1}^{-1}\ast y_{n}^{-1})^{n}]$
where $y_{i}$ loops counterclockwise around the $ith$ circle. The sequence
diverges in $\pi _{1}(X,p)$ with the quotient topology but the sequence
converges to the trivial element in the inverse limit space $%
\lim_{\leftarrow }\pi _{1}(\cup _{i=1}^{n}S_{i}).$

\section{Main Result}

Suppose $X$ is a topological space and $p\in X$. Endowed with the compact
open topology, let $C_{p}(X)=\{f:[0,1]\rightarrow X$ such that $f$ is
continuous and $f(0)=f(1)=p\}.$ Then the \textsl{topological fundamental
group} $\pi _{1}(X,p)$ is the quotient space of $C_{p}(X)$ obtained by
treating the path components of $C_{p}(X)$ as points. Thus, letting $%
q:C_{p}([0,1],X)\rightarrow \pi _{1}(X,p)$ denote the canonical surjection,
a set $A\subset \pi _{1}(X,p)$ is closed in $\pi _{1}(X,p)$ if and only if $%
q^{-1}(A)$ is closed in $C_{p}([0,1],X)$. 

The space $Y$ is said to be $T_{1}$ if the
one point subsets of $Y$ are closed.

If $A_{1},A_{2},$ are topological spaces and $f_{n}:A_{n+1}\rightarrow A_{n}$
is a continuous surjection then, (endowing $A_{1}\times A_{2}...$ with the
product topology) the \textsl{inverse limit space }$\lim_{\leftarrow
}A_{n}=\{(a_{1},a_{2},...)\in (A_{1}\times A_{2}...)|f_{n}(a_{n+1})=a_{n}\}.$

The map $f:[0,1]\rightarrow Y$ \textsl{is of the form} $\alpha _{1}\ast
\alpha _{2}...\ast \alpha _{n}$ if there exists a partition $t_{0}\leq
t_{1}...\leq t_{n}$ of $[0,1]$ such that for each $i\geq 1$ we have $%
f_{[t_{i-1},t_{i}]}=\alpha _{i}.$

For the remainder of the paper we use the following notation.

Let $X_{n}=\cup _{i=1}^{n}\{(x,y)\in R^{2}
|(x-\frac{1}{n})^{2}+y^{2}=\frac{1}{n^{2}}%
\}$. Note since $X_{n}$ is locally contractible the path components of $%
C_{p}(X_{n})$ are open in $C_{p}(X_{n})$ and hence the topological group $%
\pi _{1}(X_{n},p)$ has the discrete topology.

Let $r_{n}^{\ast }:\pi _{1}(X_{n},p)\rightarrow \pi _{1}(X_{n-1},p)$ denote
the epimorphism induced by the retraction $r_{n}:X_{n}\rightarrow X_{n-1}$
collapsing the $n^{\rm th}$ circle to the point $p=(0,0)$. Let $\lim_{\leftarrow
}\pi _{1}(X_{n},p)$ denote the inverse limit space under the maps $%
r_{n}^{\ast }.$

Let $X=\cup _{n=1}^{\infty }X_{n}$ denote the familiar Hawaiian and let $%
R_{n}:X\rightarrow X_{n}$ denote the retraction fixing $X_{n}$ pointwise and
collapsing $\cup _{i=n+1}^{\infty }X_{i}$ to the point $p.$

The formula $\phi ([f])=([R_{1}(f)],[R_{2}(f)],...)$ determines a canonical homomorphism 
$\phi :\pi 
_{1}(X,p)\rightarrow \lim_{\leftarrow }\pi _{1}(X_{n},p)$.

\theoremstyle{definition}
\newtheorem*{rem}{Remark}                % (5)
\begin{rem}
The homomorphism $\phi :\pi _{1}(X,p)\rightarrow \lim_{\leftarrow }\pi
_{1}(X_{n},p)$ is continuous (Proposition 3.3 \cite{Biss}) and one to one
(Theorem 4.1 \cite{mor}). Since $\pi _{1}(X_{n},p)$ is discrete the space $%
\Pi _{n=1}^{\infty }\pi _{1}(X_{n},p)$ is metrizable and in particular the
subspace $\lim_{\leftarrow }\pi _{1}(X_{n},p)$ is a $T_{1}$ space.
Consequently $\pi _{1}(X,p)$ is a $T_{1}$ space since $\phi $ is continuous
and one to one. Thus the path components of $C_{p}(X)$ are closed in $%
C_{p}(X).$
\end{rem}

\theoremstyle{plain}
\newtheorem{thm}{Theorem}[section]   
\begin{thm}
\label{main}The injection $\phi :\pi _{1}(X,\{p\})\hookrightarrow
\lim_{\leftarrow }\pi _{1}(X_{n},p)$ is not a topological embedding.
\end{thm}

\begin{proof}
Let $q=(2,0)$ in $X_{1}$. For a loop $f:[0,1]\rightarrow \cup
_{i=1}^{\infty }X_{i}$ with base point $p=(0,0)$ define the oscillation
number $O_{q}(f)$ as the maximal $n$ such that there exist $%
0=t_{0}<t_{1}\cdots t_{2n-1}<t_{2n}=1$ with $f(t_{2i})=p$ and $f(t_{2i-1})=q$%
. Let $y_{i}\in C_{p}(X)$ loop once counterclockwise around the $ith$ circle
and let $y_{i}^{-1}$ $\in C_{p}(X)$ loop once clockwise around the $ith$
circle.

First note that if $f\in C_{p}(\cup _{i=1}^{\infty }X_{i})$ is path
homotopic to a map of the form $(y_{1}^{-1}\ast y_{n}^{-1}\ast y_{1}\ast
y_{n})^{n}$ then $O_{q}(f)\geq 2n$ for $n\geq 2.$ To see this first observe $%
O_{q}(f)=O_{q}(R_{n}f).$ Now recall $\pi _{1}(X_{n},p)$ is canonically
isomorphic to the free group on generators $\{y_{1},...y_{n}\}.$ Thus if $w$
is an (unreduced) word corresponding to $R_{n}f$ then each step of the
algebraic reduction of $w$ to $(y_{1}^{-1}y_{n}^{-1}y_{1}y_{n})^{n}$ never
raises the oscillation number of the corresponding path in $X_{n}.$ Hence $%
O_{q}(f)\geq O_{q}((y_{1}^{-1}\ast y_{n}^{-1}\ast y_{1}\ast y_{n})^{n})=2n.$

To prove $\phi $ is not an embedding consider the set $A\subset \pi
_{1}(X,p) $ defined as $A=\{[f_{2}],[f_{3}]...]\}$ where $f_{n}$ is of the
form $(y_{1}^{-1}\ast y_{n}^{-1}\ast y_{1}\ast y_{n})^{n}.$ To see that $A$
is closed in $\pi _{1}(X,p)$ consider the union of (closed) 
path components $B=\cup
_{n=2}^{\infty }[f_{n}]\subset C_{p}(X).$ Observe if $f\in C_{p}(X)$ there exists an
open neighborhood $U\subset C_{p}(X)$ such that $O_{q}(f)\geq O_{q}(g)$ for
each $g\in U.$ Thus $U\cap \lbrack f_{n}]\neq \emptyset $ for at most
finitely many of the closed sets $[f_{n}].$ Hence $B$ is closed in $C_{p}(X)$
and consequently $A$ is closed in $\pi _{1}(X,p).$ On the other hand $\phi
(A)$ is not closed in the image of $\phi $ since the sequence $\{\phi
([f_{n}])\}$ converges to the trivial element in $\lim_{\leftarrow }\pi
_{1}(X_{n},p).$ Hence $\phi $ is not a homeomorphism from $\pi _{1}(X,p)$
onto the image of $\phi .$
\end{proof}

\Addresses\recd

\end{document}